\documentclass[12pt]{amsart}
\usepackage{latexsym, amsmath,amssymb}

\usepackage{graphicx}

\usepackage{xcolor}

 \numberwithin{equation}{section}

\renewcommand{\epsilon}{\varepsilon}
\newtheorem{theorem}{Theorem}

\newtheorem{lemma}[theorem]{Lemma}

\newtheorem{proposition}[theorem]{Proposition}

\numberwithin{theorem}{section}

\newcommand{\vertiii}[1]{{\left\vert\kern-0.25ex\left\vert\kern-0.25ex\left\vert #1 
    \right\vert\kern-0.25ex\right\vert\kern-0.25ex\right\vert}}

\begin{document}

\title[$C^{1,1-\epsilon}$ Isometric embeddings]
{$C^{1,1-\epsilon}$ Isometric embeddings}

\author[A. D. Mart\'inez]{\'Angel D. Mart\'inez}
\address{Institute for Advanced Study, Fuld Hall 412, 1 Einstein Drive, Princeton, NJ 08540, United States of America} 
\email{amartinez@ias.edu}

\maketitle

\begin{abstract}
In this paper we use the convex integration technique enhanced by an extra iteration originally due to K\"all\'en and revisited by Kr\"oner to provide a local $h$-principle for isometric embeddings in the class $C^{1,1-\epsilon}$ for $n$-dimensional manifolds in codimension $\frac{1}{2}n(n+1)$.
\end{abstract}

\section{Introduction}\label{intro}

At the beginning of the 20th century the concepts of differential geometry had been defined in two different ways: extrinsically, using some ambient space; or, intrinsically, using coordinate charts as B. Riemann did in his famous {\em \"Uber die Hypothesen welche der Geometrie zu Grunde liegen} (1854). For some time this two approaches co-existed. Until H. Whitney proved that they are equivalent in 1936. More concretely, he proved that it is nevertheless possible to smoothly embed a manifold into an Euclidean space of high dimension. This settled the equivalence of definitions for topological smooth manifolds. It was not clear, however, how to produce an embedding that preserved lengths. 

In 1827, even before B. Riemann delivered his {\em Habilitationsschrift} lecture on geometry, C. F. Gauss had proved his \textit{theorema egregium} which shows that the scalar curvature is an intrinsic object invariant by isometries. It was already known how to relate the curvature of surfaces in $\mathbb{R}^3$ to the ambient metric provided the embedding is twice differentiable and isometric. This imposes some rigidity for $C^2$ isometric maps. Indeed, a classical application of this shows the impossibility to isometrically embed a two dimensional flat torus in the three dimensional euclidean space.  It was rather natural to try to extend this rigidity to $C^1$ maps which is the minimal regularity needed to pull-back the metric from the ambient space. This remained an open problem for decades.

Until 1954, when J. F. Nash astonished the mathematical community with an ingenious construction of $C^1$ isometric embeddings for Riemannian manifolds in codimension at least two. His construction departed from a smooth embedding whose existence is assured in high codimension by Whitney's embedding theorem. (The immersion conjecture, proved by R. L. Cohen in 1985, predicts the sharp codimension that an immersion would require.) Nash also predicted that a variant of his proof could lower the codimension restriction to one. This was achieved by Kuiper in 1955. This result is nowadays known as the Nash-Kuiper theorem. Nash's proof proceeds providing an iterative scheme that converges in $C^1$ norm. Each stage consist on a series of steps where the metric tensor is succesively approximated with the help of a spiral in a plane which is essentially perpendicular to the tangent plane. Kuiper proceeded in a similar manner, replacing spirals by what he called strains  (compare with corrugations). 

A refinement of this scheme does in fact converge in $C^{1,\alpha}$ for some positive $\alpha$ which depends, for instance, on the number of steps and, therefore, deteriorates with the dimension. A $C^{1,1/7-\epsilon}$ result for analytic metrics on the disc was announced by Borisov in 1965. The best general bounds for embeddings in $\mathbb{R}^{n+1}$ of $n$-dimensional balls and smooth metrics are due to Conti, De Lellis and Sz\'ekelyhidi which give $\alpha<\frac{1}{1+n+n^2}$. One of the difficulties that arise is the loss of derivatives inherent to the method. Borisov avoided this difficulty using the analyticity hypothesis and keeping track of all the derivatives. On the other hand Conti, De Lellis and Sz\'ekelyhidi get rid of the analyticity assumption using a mollification step. 

Quite recently De Lellis, Inauen and Sz\'ekelyhidi reached $\alpha<1/5-\epsilon$ for embeddings of two dimensional balls in three dimensional space. The improvement is  possible, at least heuristically, due to the existence of conformal mappings that diagonalize the metric and, therefore, reduce the number of steps. Indeed, with the help of this change of variables, instead of approximating three tensorial  coordinates (the tensor is symmetric), one only needs to approximate the diagonal. Furthermore, the new coordinates are isothermal which implies both entries in the diagonal are equal. This technique can not be carried to higher dimensions (cf. Liouville's theorem).  

A striking conjecture due to M. Gromov predicts that the threeshold for flexibility is $\alpha<1/2$ in codimension one (cf. \cite{G}, Problem 39, also \cite{LS, DIS, Kl}). There is strong evidence in favour of this. Let us mention the work of De Lellis and Inauen were the authors found that this exponent is critical for the Levi-Civita connection and studied the problem for spherical caps in codimension twelve (cf. \cite{IL}). D. Inauen has introduced a natural extrinsic notion for parallel transport in his PhD dissertation that works in codimension one (or higher) and coincides with the intrinsic notion provided $\alpha>\frac{1}{2}(\sqrt{5}-1)>1/2$ (cf. \cite{Inauen}). This rigidity is lost in higher codimension where even better regularity can be achieved as J. F. Nash already showed (cf. \cite{N2,  Gr, Ka}). The literature is vast and we refrain from trying to provide a complete account of it here. \textcolor{black}{The best result is due to G\"unther proving the existence of smooth isometric immersions \textcolor{black}{in codimension $\frac{1}{2}n(n+1)+\max\{n,5\}$}.} We refer the reader to \cite{N, K, B, B2, Co, CLS, D, Gun, Gun2,Inauen} and the references therein. 

The following result explicitly stated in Kr\"oner's master thesis, after K\"allen's work, provides an $h$-principle in the class $C^{1,1-\epsilon}$ in a slightly smaller codimension than the best for which a smooth result is known to hold:

\begin{theorem}[K\"allen, Kr\"oner]\label{c11kroner}
Let \textcolor{black}{$\alpha<1$}, $\mu>0$  and $h:(B^n_{1+\mu},g)\rightarrow\mathbb{R}^{n+n\textcolor{black}{(n+1)}}$ be a short embedding with smooth metric $g$. Then, for any $\epsilon>0$, there is a $C^{1,\alpha}$ isometric embedding $f:(B^n_1,g)\rightarrow(\mathbb{R}^{n+\textcolor{black}{n(n+1)}}, e)$ such that 
\[\|f-h\|_{0}\leq \epsilon.\]
\end{theorem}

We will revisit the argument leading to this result which relies on an ingenious iteration introduced by K\"allen together with the vanishing of an error of a certain type. The main result in this paper improves this result by reducing the necessary codimension by half, namely \textcolor{black}{$\frac{1}{2}n(n+1)$}:

\begin{theorem}\label{first}
Let $\alpha<1$, $\mu>0$  and $h:(B^n_{1+\mu},g)\rightarrow\mathbb{R}^{n+\textcolor{black}{\frac{1}{2}n(n+1)}}$ be a short embedding with metric smooth $g$. Then, for any $\epsilon>0$, there is a $C^{1,\alpha}$ isometric embedding $f:(B^n_1,g)\rightarrow(\mathbb{R}^{n+\textcolor{black}{\frac{1}{2}n(n+1)}}, e)$ such that 
\[\|f-h\|_{0}\leq \epsilon.\]
\end{theorem}

\textcolor{black}{In this paper we apply a new variant of the strains a la Kuiper for which a miraculous cancellation occurs.} The theorem also holds for immersions although our presentation will not take care of this duplicity and will only deal with the embedding case explicitly. The proof of this result follows the standard arguments as in the literature (cf. \cite{CLS, DIS}) together with the aforementioned almost fixed point theorem due to K\"all\'en. 

It should be noted that in this paper we focus our attention to the local case only. A possible extension of the argument, following the setting of {\em adapted short immersions} introduced by Cao and Sz\'ekelyhidi, should provide the global version of this result (cf. \cite{Ka, CS3}). Another fundamental difference with K\"allen's work is that it applies for rough metrics too. We shall not attempt to achieve, nor describe this here, and will focus on the main cancellation phenomena which is new and allows to drop the dimension by a factor of two. 

For the two dimensional sphere it is known that this kind of $h$-principles can not hold if $\alpha>2/3$ in $\mathbb{R}^3$ due to a result of Borisov. Very recently De Lellis and Pakzad have showed the same limitation for isometric immersions of the flat two dimensional torus in codimension one (cf. \cite{CP}). We refer to the literature for more details on  rigidity results (cf. \cite{P} or \cite{CLS} for an alternative proof, and the references therein). 

It is quite remarkable that the convex integration technique is also useful to construct H\"older continuous solutions for the Euler equation that do not conserve energy as predicted by Onsager in 1949. This has been first observed by De Lellis and Sz\'ekelyhidi who  initiated a succesful program to attack this problem. This culminated with P. Isett's breakthrough which allowed him to construct solutions compactly supported in time with the optimal H\"older exponent $\alpha<1/3$. Recently after the publication of this result Buckmaster, De Lellis, Sz\'ekelyhidi and Vicol where able to construct solutions whose energy decays (cf. \cite{LS0, LS2, LS1, BDIS, I, BDSV}).

\section{Main ingredients and Sketch of proof}

\subsection{Heuristics of the method} \label{heuristics}

Nash's idea can be summarized as follows. Given a compact Riemannian manifold $(M,g)$ first use Whitney's theorem to produce a smooth embedding (or you might have one of your own choice \textit{a priori}). It is very unlikely that it is an isometric immersion but, by a suitable contraction of the euclidean space, one is able to produce a different embedding such that all lengths between points in the metric induced by the embedding are shorter than the ones the metric $g$ would induce. Embeddings satisfying this property are known as {\em short}. Then one runs an iterative argument that will smoothly enlongate the manifold at each point (and direction) simultaneously. Some parts will need to be enlongated more than others though and one would not expect to achieve this in a single stage. In order to make this an inductive argument one needs to end up with a short embedding again. This produces a sequence of embeddings whose induced metrics better converge to $g$. At the same time one observes that a number of derivatives of the embeddings might (and will) blow up. Nash was able to control the $C^1$ norms of this sequence in the limit. Our interest is on $C^{1,\alpha}$ norms for which some notation will be needed.

The scheme proceeds iteratively from a given short embedding $f_q:M\rightarrow \mathbb{R}^{n+1}$, such that the induced metric $g^q=f_q^{\sharp}e$ is at some distance to the desired metric, say, $\|g_{ij}-g^{q}_{ij}\|_0\leq \delta_{q+1}$. Then it provides a new approximation $f_{q+1}=f_q+w_{q+1}$, which is short again, using a spiral $w_{q+1}$, which is a vector field which oscillates at frequency $\lambda_{q+1}$ that essentially solves
\[g_{ij}-g^q_{ij}=\partial_iw_{q+1}\cdot\partial_jw_{q+1}+E_{q+1}\]
with probably some small error term $E_{q+1}$ that we can control when we choose the perturbation $w_{q+1}$. This suggests that the spiral satisfies $\|w_{q+1}\|_0\leq \delta_{q+1}^{1/2}/\lambda_{q+1}$. With this in mind let us comment on  the inductive step, which reduces to bound the difference
\[g_{ij}-g_{ij}^{q+1}=g_{ij}-g_{ij}^q-2\partial_if_q\cdot\partial_jw_{q+1}-\partial_iw_{q+1}\cdot\partial_jw_{q+1}\]
by $\delta_{q+2}$. It would be enough for $\partial_iw_{q+1}$ to be essentially orthogonal to the image of $f_q$ and the error term $E_{q+1}$ in the previous equation to be smaller than $\delta_{q+2}$. A technical requirement to close the argument is to ensure $f_{q+1}$ is again a short embedding.

It is worth mentioning now how to control the $C^{1,\alpha'}$ regularity of the limit, $f$. For the sake of illustration let us employ the \textit{geometric ansatz}, namely $\delta_q=\lambda_q^{-2\alpha}=\lambda^{-2\alpha q}$ for some $\lambda>1$, until the end of this section. Notice that $\delta_q$ tends to zero which would imply that the metric at each stage $g^q$ tends uniformly to the desired metric $g$. On the other hand it is clear that
\[\|f-f_0\|_{C^{1,\alpha'}}\leq\sum_{q=1}^{\infty}\|w_q\|_{1,\alpha'}=\sum_{q=1}^{\infty}\delta_q^{1/2}\lambda_q^{\alpha}=\sum_{q=1}^{\infty}\lambda^{(\alpha'-\alpha)q}\]
which converges if $\alpha'<\alpha$. 

Let us finally state a variant of the H\"older version of the Nash-Kuiper theorem for a ball in codimension one:

\begin{theorem}[cf. \cite{CLS}]\label{nashkuiper}
Let $\eta\geq 0$, $\mu>0$ and $h:(B^n_{1+\mu},g)\rightarrow(\mathbb{R}^{n+1},e)$ be a short embedding such that $g-h^{\sharp}e\geq\eta\textrm{Id}$. Then there is some $\alpha_0(n)>0$ such that for any $\epsilon>0$ there exist a $C^{1,\alpha_0}$ embedding $h_0:(B^n_1,g)\rightarrow\mathbb{R}^{n+1}$ such that $\|h-h_0\|_{0}\leq\epsilon$ and $g-h_0^{\sharp}e=\eta\textrm{Id}$.
\end{theorem}

Notice that this corresponds to the work of Conti, De Lellis and Szek\'elyhidi for $\eta=0$. In fact it is equivalent to it. Indeed, notice that $\tilde{g}=g-\eta\textrm{Id}\geq h^{\sharp}e\geq 0$ is a metric for which $h$ is a short embedding. Then there is an isometric embedding $h_0:(B,\tilde{g})\rightarrow(\mathbb{R}^{n+1},e)$ which satisfies the conclusion above. \textcolor{black}{We refer the reader to the literature for a proof.} We include this statement here for reference only as it will be useful later to provide an initial step for the induction process that proves Theorem \ref{first} (cf. Section \ref{prooffirst}).

In the next section we provide a statement of an inductive step that is key in the proof. In Section \ref{lemmata} we resume some technical auxiliary results that will be used in Section \ref{proofproposition}. Section \ref{theperturbation} is devoted to introduce the general form of the perturbation and some useful notation . In Section \ref{prooffirst} we conclude the proof of Theorem \ref{c11kroner}. Finally, in Section \ref{theperturbationone}, we introduce the strains and explain the modifications in the argument that allow to prove Theorem \ref{first}.

\section{The inductive step}\label{inductivestep}

The inductive step will need a specific choice of parameters to work. In this section we present the so-called \textit{double exponential ansatz}. Namely,
\[\textrm{$\delta_{q}=a^{-2\alpha b^{q-1}}$ and $\lambda_q=a^{b^q}$}.\]
(This is needed due to the $\epsilon$ loss in our choice of $\ell_q$, defined in Section \ref{lemmata}). The parameter $a$ will be chosen large enough and $b=b(\alpha)>1$ independently of the stage $q$. Notice that under this ansatz
\[\|f-f_0\|_{C^{1,\alpha'}}\leq\sum_{q=1}^{\infty}\|w_q\|_{1,\alpha'}=\sum_{q=1}^{\infty}\delta_q^{1/2}\lambda_q^{\alpha}=\sum_{q=1}^{\infty}a^{-(\alpha-\alpha' b)b^{q-1}}\]
which converges if $\alpha'<\alpha/b$. This already suggests the need to choose $b$ close to one to achieve optimal results.

The following result provides an inductive step that allows to control the H\"older norms through the iteration

\begin{proposition}\label{induction}
Let $g$ be the smooth fixed metric, \textcolor{black}{$\beta\in(0,1)$ and 
\begin{equation}\label{restriction2}
2\alpha<2-\beta
\end{equation}
For a sufficiently small $\epsilon\in(0,1/4)$ there exist a choice of parameters in the double exponential ansatz such that the following holds.} Suppose that we have constructed smooth short embeddings $f_0,\ldots,f_q:(B_R,g)\rightarrow(\mathbb{R}^{n+\frac{n(n+1)}{2}}, e)$ such that $f_{q}=f^{\ell_q}_{q-1}+w_q$ (cf. Section \ref{proofproposition}) where
\[\left\|\frac{g_{ij}-g_{ij}^q}{\delta_{q+1}}-\operatorname{Id}\right\|_0\leq \lambda_q^{-\epsilon},\]
\[\left\|\frac{g_{ij}-g_{ij}^q}{\delta_{q+1}}-\operatorname{Id}\right\|_{\beta}\leq \lambda_q^{\beta-\epsilon},\]
\[\|w_{q}\|_r\leq C\delta_q^{1/2}\lambda_q^{r-1},\]
\[\|f_q\|_{r}\leq \kappa^{r}_q,\]
and $w_0=0$. Then there exists a smooth short embedding $f_{q+1}:(B_{R-\ell_q},g)\rightarrow(\mathbb{R}^{n+\frac{n(n+1)}{2}},e)$ such that 
\[\left\|\frac{g_{ij}-g_{ij}^{q+1}}{\delta_{q+2}}-\operatorname{Id}\right\|_0\leq \lambda_{q+1}^{-\epsilon},\]
\[\left\|\frac{g_{ij}-g_{ij}^{q+1}}{\delta_{q+2}}-\operatorname{Id}\right\|_{\beta}\leq \lambda_{q+1}^{\beta-\epsilon},\]
\[\|w_{q+1}\|_r\leq C\delta_{q+1}^{1/2}\lambda_{q+1}^{r-1}\]
and
\[\|f_{q+1}\|_{r}\leq \kappa^{r}_q+C\delta_{q+1}^{1/2}\lambda_{q+1}^{r-1}\]
for any $r\in [0,2]$ and a constant $C>1$ independent of $q$.
\end{proposition}

The hypothesis on the metric ensure the existence of a decomposition satisfying Lemma \ref{decomposition}. It is worth mentioning that
\[\left\|\frac{g_{ij}-g_{ij}^q}{\delta_{q+1}}-\textrm{Id}\right\|_0\leq \lambda_q^{-\epsilon}\]
implies
\[\|g_{ij}-g_{ij}^q\|_0\leq 2\delta_{q+1}\]
which relates to the heuristics as discussed in Section \ref{heuristics}. Furthermore, it immediately implies $f_{q+1}$ is short (cf. Lemma \ref{h}). 

 The parameter $b>1$ will be chosen sufficiently close to one, depending on $\alpha$, while $\epsilon$ will be fixed and sufficiently small depending on it. The statement implicitly assumes that $a$ is chosen large enough depending on the rest of parameters, the initial data $f_0$ and $g$. In the exposition we will deduce that such a choice is possible.  

\section{Technical lemmata}\label{lemmata}

Some remarks are in order. The proof has been partitioned in a series of lemmata below. Implicitly all of them assume that we can choose $a$ large enough and $b>1$ so that the parameters satify certain natural restrictions steming from the proofs. Let us also remark that, as the reader might expect from the statement, we will use a complete induction argument (cf. Lemma \ref{basis}). That means that we will assume the hypothesis for $f_0$,\ldots, $f_q$ and then prove it for $f_{q+1}$ (cf. Lemma \ref{basis} below). No constant in this section depends on $q$ unless we indicate it explicitly. 

Furthermore, the notation of this section will be kept to a minimum whenever possible. As it is customary the constants might change from line to line. We will therefore commit some abuse of notation for the {\em sake of clarity}. Some quantities are tensors but there are many parameters which are not tensorial such as the mollification parameter $\ell$, which will oscillate from top to bottom in the notation. In fact, in this section, and until the end of the proof of Proposition \ref{induction} in Section \ref{proofproposition} we fix
\[\ell_q^{2-\beta}=\lambda_q^{\beta-2-\epsilon}\delta_{q+1}\delta_q^{-1}\]
but we will omit the dependence on $q$ (which is fixed, cf. Proposition \ref{induction}) and denote it as $\ell=\ell_q$. Another simplification is that we will not keep track on the shrinking process but will indicate when it occurs for the reader's convenience. This is not dangerous because the total shrinking length, namely
\[\sum_{q=0}^{\infty}\ell_q\leq\sum_{q=0}^{\infty}\lambda_q^{-1},\]
can be made smaller than $\mu$ by choosing $a$ large enough (cf. Theorem \ref{first}). Analogously, $\kappa^1_q$ is bounded uniformly by a constant depending on the initial $f_0$ only and $\kappa^2_q\leq C\delta_q^{1/2}\lambda_q$ for some appropiate constant $C>0$ independent of $q$. Finally, at many points the exposition is purposely redundant for the very same reason. We hope all this does not affect readability but enhance it. 

We denote the mollification of $f_{q}$ at scale $\ell$ by $f_q^{\ell}$ and the metric it induces, $(f^{\ell}_q)^{\sharp}e$, by $g_{\ell}^{q}$ (cf. Section \ref{sectionholder}). The following clarifies our choice for $\ell$.

\begin{lemma}[Mollification]\label{mollification}
Under the hypothesis of Proposition \ref{induction}. The mollified metric satisfies
\[\left\|\frac{g-g^q_{\ell}}{\delta_{q+1}}-\operatorname{Id}\right\|_0\leq C\lambda_q^{-\epsilon}(\lambda_q\ell)^{\beta}\leq C\lambda_q^{-\epsilon} \]
and, more generally,  
\[\left\|\frac{g-g^q_{\ell}}{\delta_{q+1}}-\operatorname{Id}\right\|_r\leq C\lambda_q^{\beta-\epsilon}\ell^{-r+\beta}\]
for any $r\geq 0$. We can also control the mollification as
\[\|f_q^{\ell}\|_{r}\leq \kappa_q^r\textrm{ for $r=1,2$}\]
and
\[\|f_q^{\ell}\|_r\leq C\delta_q^{1/2}\lambda_q\ell^{-r+2}\textrm{ for $r\geq 2$}\]
where the constants depend only on the constants from Lemma \ref{holder}, $r$ and $\|g\|_{r+2}$.
\end{lemma}

\textsc{Proof:} observe that once it is proved for $r\in\mathbb{N}$ the rest follows by interpolation (cf. Section \ref{sectionholder}). By definition, 
\[\|g-g_{\ell}^q-\delta_{q+1}\operatorname{Id}\|_r=\|g_{ij}-\partial_if_{q}^{\ell}\cdot\partial_jf_{q}^{\ell}-\delta_{q+1}\textrm{Id}_{ij}\|_r\]
which can be splitted into  three terms as follows
\[\|[g_{ij}-g_{ij}*\varphi_{\ell}]+[g_{ij}-(\partial_if_{q}\cdot\partial_jf_{q})-\delta_{q+1}\textrm{Id}_{ij}]*\varphi_{\ell}+[(\partial_if_{q}\cdot\partial_jf_{q})*\varphi_{\ell}-\partial_if_{q}^{\ell}\cdot \partial_jf_{q}^{\ell}]\|_r.\]
This is where the domain has to be shrinked as in the statement of Proposition \ref{induction} in order for the convolution to make sense. Using the bounds from Lemma \ref{holder}, this can be bounded from above by
\begin{equation}\label{moleq}
 C\left(\ell^2\|g_{ij}\|_{r+2}+\ell^{-r+\beta}\|g_{ij}-\partial_if_{q}\cdot\partial_jf_{q}-\delta_{q+1}\textrm{Id}\|_{\beta}+\ell^{2\alpha-r}\|f_{q}\|_{1,\alpha}^2\right)
\end{equation}
for any $\alpha\in(0,1]$. Choosing $\alpha=1$ the above is bounded by
\[ C(\ell^2\|g\|_{r+2}+\delta_{q+1}\ell^{\beta-r}\lambda_q^{\beta-\epsilon}+\ell^{2-r}\delta_{q}\lambda_{q}^{2})\leq C\lambda_{q}^{\beta-\epsilon}\ell^{\beta-r}\delta_{q+1}.\]
Indeed, absorbing the constant $\|g\|_{r+2}$ in $C$, the first term is smaller than the right hand side because  $\ell^2\leq \delta_{q+1}\lambda_q^{\beta-\epsilon}\ell^{\beta}$ holds or, equivalently, after some algebraic manipulation,
\[\delta_q^{-1}\leq\lambda_q^2.\]
Assuming the double exponential ansatz this follows if  $2\alpha<2b$, which holds for any $b>1$. On the other hand, the second and third terms balance equality for our particular choice of $\ell$. The last two inequalities follow from Lemma \ref{holder} (i) inmediately. For instance,
\[\|f^{\ell}_q\|_r=\|f_q*\varphi_{\ell}\|_r\leq C\ell^{-r+2} \|f_q\|_2.\]
 It will be useful to denote 
\begin{equation}\label{ell}
\ell=\lambda_q^{-1-\epsilon^*}\left(\frac{\delta_{q+1}}{\delta_q}\right)^{\beta^*}
\end{equation}
where $\beta^*=(2-\beta)^{-1}$ and $\epsilon^*=\epsilon\beta^*$ are positive constants. From this it is evident that $\lambda_q\ell\leq 1$ which shows that we are really perturbing our hypothesis in Proposition \ref{induction} at the same scale,  for $r=0,\beta$ respectively, finishing the proof.

In this case we shall use a metric decomposition due to K\"all\'en instead of isothermal coordinates (cf. \cite{paper3.0}).

\begin{lemma}[Metric decomposition]\label{decomposition}

For any fixed vectors $n_0,\ldots,\textcolor{black}{n_{n(n+1)/2-1}}\in\mathbb{S}^{n-1}$, there exists a $\sigma_1>0$ with the following property: given  \textcolor{black}{symmetric tensor} valued functions $\tau$, $\tau_k$  and $\tau_{kk'}$  such that 
\[\|\tau-\operatorname{Id}\|_0+\sum_{k=0}^{\frac{n(n+1)}{2}-1}\|\tau_{k}\|_0+\sum_{k,k'=0}^{\frac{n(n+1)}{2}-1}\|\tau_{kk'}\|_0\leq\sigma_1,\]
 then there exist a vector valued function $A$ such that
\[\tau(x)=\sum_{k=0}^{\frac{n(n+1)}{2}-1}A_k(x)^2 n_k\otimes n_k+\sum_{k=0}^{\frac{n(n+1)}{2}-1}A_k(x)\tau_{k}(x)+\sum_{k,k'=0}^{\frac{n(n+1)}{2}-1}A_k(x)A_{k'}(x)\tau_{kk'}(x)\]
with $c_i>\sigma_1$ and satisfying the estimate
\[\|A_i\|_r\leq C_r\left(1+\|\tau\|_r+\sum_{k=0}^{\frac{n(n+1)}{2}-1}\|\tau_{k}\|_r+\sum_{k,k'=0}^{\frac{n(n+1)}{2}-1}\|\tau_{kk'}\|_r\right)\] 
for any $r\in\mathbb{N}$. Furthermore, there exist a smooth $F$ such that 
\[A_k(x)=F(\tau(x), \tau_k(x), \tau_{kk'}(x)).\]
\end{lemma}

\textcolor{black}{A recent proof of this fact, using the implicit function theorem, can be found in \cite{IL} (cf. Proposition 5.4). Their statement does not correspond to the one we provide here but their proof does provide the extra information on the existence of a smooth $F$ from which it is derived.}

\textcolor{black}{In the next result we introduce bounds for the metric error $h$.}

\begin{lemma}\label{h}
Under the hypothesis of Proposition \ref{induction}. For any $\theta\in(0,1)$ the tensor $h$ defined by $\delta_{q+1}h=g-g^q_{\ell}-\delta_{q+2}\operatorname{Id}$ is positive and satisfies
\[\|h-\operatorname{Id}\|_{\epsilon^2}<\theta.\]
Furthermore, 
\[\|h\|_r\leq C\lambda_q^{\beta-\epsilon}\ell^{-r+\beta}\]
for any $r\in \mathbb{N}$.
\end{lemma}

\textsc{Remark:} $\theta$ is a parameter that will be chosen later small enough and independently of $q$. 

\textsc{Proof:} indeed, if it is not positive at some point there will be a unit vector $v$ such that the quadratic form is non positive, i.e. 
\[0\geq v^Thv=v^T[\delta_{q+1}^{-1}(g-g^q_{\ell})-\operatorname{Id}]v+(1-\delta_{q+2}\delta_{q+1}^{-1})v^T\textrm{Id} v.\]
Lemma \ref{mollification} implies that the first is bounded by $O(\lambda_q^{-\epsilon})$ while the second is close to one since $\delta_{q+2}$ decays much faster than $\delta_{q+1}$. Indeed, $\delta_{q+2}\delta_{q+1}^{-1}=o(1)$ if $-2\alpha(b-1)<0$ which follows from $b>1$. Here we are using the notation $o(1)$ to denote a quantity that tends to zero as $a$ increases to infinity. This implies $0\geq 1+o(1)$, a contradiction. On the other hand let us observe that Lemma \ref{mollification} implies
\[\|h-\textrm{Id}\|_{\epsilon^2}\leq\|(\delta_{q+1}^{-1}(g-g^q_{\ell})-\textrm{Id})-\delta_{q+2}\delta_{q+1}^{-1}\textrm{Id}\|_{\epsilon^2}\leq C\ell^{-\epsilon^2}\lambda_{q}^{-\epsilon}+\delta_{q+2}\delta_{q+1}^{-1}<\theta\]
which is true if $\lambda_{q}^{-\epsilon}\ell^{-\epsilon^2}=o(1)$. Using this the notation above this is equivalent to
\[\lambda_{q}^{-1}(\lambda_{q}^{1+\epsilon^*}\delta_{q+1}^{-\beta^*}\delta_q^{\beta^*})^{\epsilon}=o(1).\]
 Using the double exponential ansatz this reads
\[-b+\epsilon((1+2\alpha\beta^*+\epsilon^*)b-2\alpha\beta^*)< 0\]
which holds for any $\alpha\in (0,1)$ and $b>1$ provided $\epsilon<1/4$ is small enough. This forces $a=a(\epsilon)$ to be large. The last part follows by a direct application of the definition and the estimates from Lemma \ref{mollification}.

We can take $\theta<\sigma_1/2$ as in Lemma \ref{decomposition}. This shows that $h$, defined as in Lemma \ref{h} above, can be considered a metric satisfying the hypothesis of Lemma \ref{decomposition}. 

To construct the spirals we will need a pair of smooth vector fields orthogonal to the surface at any point. The following result provides their existence together with some H\"older estimates we will need later.

\begin{lemma}\label{basis}
Let $f_0$ be a smooth immersion \textcolor{black}{into $\mathbb{R}^{n+d}$.} There exist a constant $\sigma_0$ such that for any smooth map $f$ such that $\|f-f_0\|_1<\sigma_0$. Then there exist smooth orthonormal vector fields $\nu^f_i(x)$ for $i=1\textcolor{black}{,\ldots, d}$ such that
\begin{itemize}
\item[(i)] They are orthogonal to the plane tangent to the image of $f$ at $x$. 
\item[(ii)] Furthermore, they satisfy the estimate \[\|\nu^f\|_r\leq C_r(1+\|f\|_{r+1}).\]
\end{itemize}
\end{lemma}
\textsc{Remark:} this corresponds to Proposition 5.3 from \cite{IL} where the reader will find a detailed proof. The constants involved, including $\sigma_0$, depend on a finite number of derivatives of $f_0$. 

\section{The perturbation $w_{q+1}$}\label{theperturbation}

The perturbation $w_{q+1}$ will be defined inductively. First we define $f^{(0)}=f_q^{\ell}$ and $f^{(k)}=f^{(k-1)}+w_{q+1}^{(k-1)}$ for $k\geq 1$. \textcolor{black}{To define the perturbation $w^{(k)}_{q+1}$ we will need an orthonormal basis $\nu^i_k$, $i=1,2$, \textcolor{black}{$k=0,\ldots,\frac{n(n+1)}{2}-1$} orthogonal to the tangent plane of $f^{(0)}$, as given by Lemma \ref{basis} with $d=n(n+1)$. Then we can define}
\[w^{(k)}_{q+1}(x)=\delta^{1/2}_{q+1}\frac{A_k(x)}{\lambda_{q+1}}\left(\sin(\lambda_{q+1}x\cdot n_{k})\nu^1_{k}(x)+\cos(\lambda_{q+1}x\cdot n_{k})\nu^2_{k}(x)\right)\]
where the frequencies are given by \textcolor{black}{$\lambda_{q+1}=\lambda_q^{b}$.}  The choice of coefficients $A_k$ will be specified later (in the proof of Lemma \ref{perturbation} below). For the time being let us advance that  the perturbation will have the form
\[w_{q+1}(x)=\delta^{1/2}_{q+1}\sum_{k=0}^{\frac{n(n+1)}{2}-1}\frac{A_k(x)}{\lambda_{q+1}}\left(\sin(\lambda_{q+1}x\cdot n_{k})\nu^1_{k}(x)+\cos(\lambda_{q+1}x\cdot n_{k})\nu^2_{k}(x)\right).\]
Before proceeding with the proof itself we we need to introduce some notation. Let us observe that  by construction
\[\partial_i f_{q+1}\cdot\partial_j f_{q+1}=\partial_i f^{\ell}_{q}\cdot\partial_j f^{\ell}_{q}+\sum_{k=0}^{\frac{n(n+1)-1}{2}}\partial_i f_q^{\ell}\cdot\partial_j w^{(k)}_{q+1}+\sum_{k,k'=0}^{\frac{n(n+1)-1}{2}}\partial_i w_{q+1}^{(k)}\cdot\partial_j w^{(k')}_{q+1}\]
holds. Expanding the perturbations one observes the identity
\[\begin{aligned}
g^{q+1}_{ij}&=(g_{\ell}^q)_{ij}+\sum_{k=0}^{\frac{n(n+1)}{2}-1}\partial_if_q^{\ell}\cdot\partial_j\left(\delta_{q+1}^{1/2}\frac{A_k(x)}{\lambda_{q+1}}\sin(\lambda_{q+1}x\cdot n_{k})\right)\nu^1_{k}(x)\\
&\qquad\qquad\qquad\qquad\quad+\partial_i\textcolor{black}{f_q^{\ell}}\cdot\partial_j\left(\delta_{q+1}^{1/2}\frac{A_k(x)}{\lambda_{q+1}}\cos(\lambda_{q+1}x\cdot n_{k})\right)\nu^2_{k}(x)\\
&+\sum_{k=0}^{\frac{n(n+1)}{2}-1}\partial_i\textcolor{black}{f_q^{\ell}}\cdot\delta_{q+1}^{1/2}\frac{A_k(x)}{\lambda_{q+1}}\left(\partial_j\nu^1_{k}(x)\sin(\lambda_{q+1}x\cdot n_{k})+\partial_j\nu^2_{k}(x)\cos(\lambda_{q+1}x\cdot n_{k})\right)\\
&+\sum_{k\textcolor{black}{,k'}=0}^{\frac{n(n+1)}{2}-1}\partial_iw^{(k)}_{q+1}\cdot\partial_jw^{(\textcolor{black}{k'})}_{q+1}
\end{aligned}\]
Notice that the first sum vanishes because the normals are orthogonal to the tangent to $f_q^{\ell}$ at $x$. We are abusing the notation, as we will, since one more term, corresponding to the second summand with $i$ and $j$ reversed, should appear (cf. Section \ref{theperturbationone} below). It will be convenient to introduce the following notation for a relevant linear term (which is part of the second sum above)
\[L(A)=\sum_{k=0}^{\frac{n(n+1)}{2}-1}\textcolor{black}{\partial_if_{q}^{\ell}}\cdot\frac{A_k(x)}{\lambda_{q+1}\delta_{q+1}^{1/2}}\left(\partial_j\nu^1_{k}(x)\sin(\lambda_{q+1}x\cdot n_{k})+\partial_j\nu^2_{k}(x)\cos(\lambda_{q+1}x\cdot n_{k})\right)\]
Then we get that
\[g^{q+1}=g^q+\delta_{q+1}L+\sum_{k\textcolor{black}{,k'}=0}^{\frac{n(n+1)}{2}-1}\partial_iw^{(k)}_{q+1}\cdot\partial_jw^{\textcolor{black}{(k')}}_{q+1},\]
using the notation we have just introduced for the linear term. We shall need some more notation to decompose the quadratic term. Let us digress and introduce it. The main term has the form
\begin{equation}\label{mainterm}
\begin{aligned}
M(a)=\sum_{k,k'=0}^{\frac{n(n+1)}{2}-1}&a_k(x)\left(\cos(\lambda_{q+1}x\cdot n_{k})\nu^1_{k}(x)- \sin(\lambda_{q+1}x\cdot n_{k})\nu^2_{k}(x) \right) \cdot\\
 & \qquad\cdot a_{k'}(x)\left(\cos(\lambda_{q+1}x\cdot n_{k'})\nu^1_{k'}(x)-\sin(\lambda_{q+1}x\cdot n_{k'})\nu^2_{k'}(x)\right)n_k\otimes n_{k'}\\
&=\sum_{k=0}^{\frac{n(n+1)}{2}-1}a_k(x)^2n_k\otimes n_{\textcolor{black}{k}}
\end{aligned}
\end{equation}
\textcolor{black}{using the orthogonality of the normals}. It will be convenient to introduce also some notation for a relevant bilinear tensor valued error term
$R=R_2+R_3+R_4$ where
\[\begin{aligned}
R_2(a)_{ij}=\sum_{k,k'=0}^{\frac{n(n+1)}{2}-1}&\frac{a_k(x)}{\lambda_{q+1}}\frac{a_{k'}(x)}{\lambda_{q+1}}\left(\sin(\lambda_{q+1}x\cdot n_{k})\partial_i\nu^1_{k}(x)+  \cos(\lambda_{q+1}x\cdot n_{k})\partial_i\nu^2_{k}(x)\right) \cdot\\
 & \qquad\cdot\left(\sin(\lambda_{q+1}x\cdot n_{k'})\partial_j\nu^1_{k'}(x)+\cos(\lambda_{q+1}x\cdot n_{k'})\partial_j\nu^2_{k'}(x)\right),
\end{aligned}\]
\[\begin{aligned}
R_3(a)_{ij}=\sum_{k,k'=0}^{\frac{n(n+1)}{2}-1}&\frac{\partial_ia_k(x)\partial_ja_{k'}(x)}{\lambda_{q+1}\lambda_{q+1}}\left(\sin(\lambda_{q+1}x\cdot n_{k})\nu^1_{k}(x)+\cos(\lambda_{q+1}x\cdot n_{k})\nu^2_{k}(x) \right) \cdot\\
 & \qquad\cdot\left(\sin(\lambda_{q+1}x\cdot n_{k'})\nu^1_{k'}(x)+\cos(\lambda_{q+1}x\cdot n_{k'})\nu^2_{k'}(x)\right),
\end{aligned}\]
and 
\[\begin{aligned}
R_4(a)_{ij}=\sum_{k,k'=0}^{\frac{n(n+1)}{2}-1}&\frac{a_k(x)}{\lambda_{q+1}}\frac{\partial_ja_{k'}(x)}{\lambda_{q+1}}\left(\sin(\lambda_{q+1}x\cdot n_{k})\partial_i\nu^1_{k}(x)+  \cos(\lambda_{q+1}x\cdot n_{k})\partial_i\nu^2_{k}(x)\right) \cdot\\
 & \qquad\cdot\left( \sin(\lambda_{q+1}x\cdot n_{k'})\nu^1_{k'}(x)+\cos(\lambda_{q+1}x\cdot n_{k'})\nu^2_{k'}(x)\right)
\end{aligned}.\]
Furthermore, we define
\[\begin{aligned}
R_5(a)_{ij}=\sum_{k,k'=0}^{\frac{n(n+1)}{2}-1}&\frac{a_k(x)a_{k'}(x)}{\lambda_{q+1}}\left(\sin(\lambda_{q+1}x\cdot n_{k})\partial_i\nu^1_{k}(x)+  \cos(\lambda_{q+1}x\cdot n_{k})\partial_i\nu^2_{k}(x)\right) \cdot\\
 & \qquad\cdot\left(\cos(\lambda_{q+1}x\cdot n_{k'})\nu^1_{k'}(x)-\sin(\lambda_{q+1}x\cdot n_{k'})\nu^2_{k'}(x)\right)(n_{k'})_j
\end{aligned}\]
and finally
\[\begin{aligned}
R_6(a)_{ij}&=\sum_{\textcolor{black}{k,k'=0}}^{\frac{n(n+1)}{2}-1}\frac{\partial_ia_k(x)a_{k'}(x)}{\lambda_{q+1}}\left(\sin(\lambda_{q+1}x\cdot n_{k})\nu^1_{k}(x)+ \cos(\lambda_{q+1}x\cdot n_{k})\nu^2_{k}(x) \right) \cdot\\
 & \qquad\qquad\qquad\cdot\left(\cos(\lambda_{q+1}x\cdot n_{k'})\nu^1_{k'}(x)-\sin(\lambda_{q+1}x\cdot n_{k'})\nu^2_{k'}(x)\right)(n_{k'})_j
\end{aligned}\]
Quite strikingly $R_6=0$ by orthogonality considerations. The main purpose of  \textcolor{black}{this analysis is to express the  quadratic term in the expression for} $g^{q+1}$ above as follows
\[\delta_{q+1}^{-1}\sum_{k\textcolor{black}{,k'}=0}^{\frac{n(n+1)}{2}-1}\partial_iw^{(k)}_{q+1}\cdot\partial_jw^{\textcolor{black}{(k')}}_{q+1}=M(A)+R_2(A)+R_3(A)+R_4(A)+R_5(A)
\]
Let us digress from the presentation a bit to make a number of remarks. The last term $R_6$ is of utmost importance in the analysis of this problem. For example, the key argument of K\"all\'en  takes advantage of the extra cancellation of $R_6$ that occurs in codimension $n(n+1)$ to prove a $C^{1,1-\epsilon}$ isometric embedding theorem (cf. \cite{Kroner}). \textcolor{black}{We shall provide a proof of Theorem \ref{c11kroner} isometric embedding theorem first and then we will show how to adapt it to reduce the codimension by half using strains.}

\section{Proof of Proposition \ref{induction}}\label{proofproposition}

In this section we construct the coefficients that we shall use to construct the needed perturbation as in our previous section. This will require the following
\begin{lemma}\label{tau}
\textcolor{black}{The tensors
\[(\tau_k)_{ij}=\frac{\partial_if_q^{\ell}}{\lambda_{q+1}\delta_{q+1}^{1/2}}\cdot\left(\partial_j\nu^{1}_{k}(x)\sin(\lambda_{q+1}x\cdot n_{k})+\partial_j\nu^{2}_{k}(x)\cos(\lambda_{q+1}x\cdot n_{k})\right)\]
and 
\[\begin{aligned}
(\tau_{kk'})_{ij}=&\frac{1}{\lambda_{q+1}}\left(\sin(\lambda_{q+1}x\cdot n_{k})\partial_i\nu^{1}_{k}(x)+  \cos(\lambda_{q+1}x\cdot n_{k})\partial_i\nu^{2}_{k}(x)\right) \cdot\\
 & \qquad\cdot\left(\cos(\lambda_{q+1}x\cdot n_{k'})\nu^{1}_{k'}(x)-\sin(\lambda_{q+1}x\cdot n_{k'})\nu^{2}_{k'}(x)\right)(n_{k'})_j \end{aligned}\]
 satisfy the estimates
 \[\|\tau_k\|_r+\|\tau_{kk'}\|_r\leq\frac{\lambda_{q+1}^r}{\lambda_{q+1}\ell}.\]}
\end{lemma}

\textsc{Proof:} it is elementary to check \textcolor{black}{using Lemmata \ref{mollification}, \ref{basis} and $\ell^{-1}\leq\lambda_{q+1}$ (cf. equation \ref{bound2} below) that}
\[\|\tau_k\|_r=\frac{O(1)}{\lambda_{q+1}\delta_{q+1}^{1/2}}\delta_q^{1/2}\lambda_q\lambda_{q+1}^r\]
while
\[\|\tau_{kk'}\|_r=\frac{O(1)}{\lambda_{q+1}}\delta_q^{1/2}\lambda_q\lambda_{q+1}^r.\]
Then the statement follows provided 
\begin{equation}\label{star}
\delta_q^{1/2}\lambda_q\leq \delta_q^{1/2}\lambda_q\delta_{q+1}^{-1/2}\leq\ell^{-1}
\end{equation}
holds. Using the definition of $\ell$ (cf. equation \ref{ell}) one can check that this is a consequence of the stronger
\[\left(\frac{\delta_{q+1}}{\delta_q}\right)^{\beta^*}\leq \left(\frac{\delta_{q+1}}{\delta_q}\right)^{1/2}\]
which holds observing that $\delta_{q+1}\leq\delta_q$ and $\beta^*\in[\frac{1}{2},1)$ by definition.

The next result provides the usual bounds for a first guess of coefficients. These approximation only takes into account $M$, $L$ and $R_5$. This will feed the ingenious iteration originally due to K\"all\'en, and revisited by Kr\"oner, that will provide the actual choice of the coefficients that we need (cf. Lemma \ref{coefficients} below).

\begin{lemma}[First coefficients]\label{conformal}
Consider Lemma \ref{decomposition} with $\tau_k$ and $\tau_{kk'}$ as in Lemma \ref{tau} and $\tau=h$ as in Lemma \ref{h}. There exist constants $C>0$ such that the following bounds hold
\begin{itemize}
\item[(i)] $\|A^{(1)}_k\|_0\leq C$.
\item[(ii)] For any $1\leq r\leq r_0$
\[\|A^{(1)}_k \|_r\leq \textcolor{black}{C\frac{\lambda_{q+1}^r}{\lambda_{q+1}\ell}}\]
\end{itemize}
\end{lemma}

\textsc{Remark:} the constants depend on $\epsilon$, $\theta$, $r_0$, the initial data $f_0$ but not on  $q$. 

\textsc{Proof:} the inequalities regarding $A$ follow from Lemmata \ref{decomposition}, \ref{tau} and \ref{h} immediately. Indeed, the bounds from Lemma \ref{tau} dominate both the constant and the ones coming from the metric error $h$. We leave further details to the reader. 

Notice that this choice of coefficients is far from enough for our purposes since the estimates do not close for the rest of error terms. To deal with this issue, improving the error, we shall use the K\"all\'en's iteration. 

\begin{proposition}[K\"all\'en's iteration]\label{kroner}
Fix $r_1\in\mathbb{N}$, $i_0$, $b$ (with inverse $F$)  and $R=R_2+R_3+R_4$ a combination bilinear terms. Suppose that for two fixed constants $\lambda,L>0$  the first $r_0+1$ derivatives of $F$ satisfy
\[\|F(T)\|_r\leq C_*\left(\|T\|_r+\frac{\lambda^r}{\lambda L}\right)\]
and
\[\|F(T)-F(T')\|_r\leq C_{*}\left(\|T-T'\|_r+(\|T'\|_r+\lambda^r)\|T-T'\|_0\right)\]
for any tensors $T$ and $T'$ in in a $(3C_*)^{-1}$-neighbourhood of $T_0$ and, futhermore, $\lambda L$ is large enough (depeding on $C_*$). Consider also any tensor $T$ in the $(3C_*)^{-1}$-neighbourhood of $T_0$ satisfying 
\[\|T\|_r\leq C\frac{\lambda^r}{\lambda L}\]
for $r\geq 1$. There exists $r_0=r_0(i_0,r_1)\in\mathbb{N}$ such that if $a^{(0)}=0$ and $a^{(1)}=F^1(T)$
satisfy
\begin{itemize}
\item[(i)] $\|a^{(1)}\|_0\leq C$.
\item[(ii)] For any $1\leq r\leq r_0-1$ the norm 
\[\|a^{(1)}\|_r\leq C\frac{\lambda^r}{\lambda L}\]
\item[(iii)] while the error
\[E_s=T-b(a^{(s)},a^{(s)})-R(a^{(s)})\]
 satisfies
\[\|E_1\|_r\leq C\frac{\lambda^r}{\lambda L}.\]
\end{itemize}
Then, there exist $a^{(i_0)}$ satisfying
\begin{itemize}
\item[(1)] $\|a^{(i_0)}\|_0\leq C$.
\item[(2)] For any $1\leq r\leq r_1$ the norm 
\[\|a^{(i_0)}\|_r\leq C'\frac{\lambda^r}{\lambda L}\]
\item[(3)] while the error satisfies
\[\|E_{i_0}\|_r\leq C'\frac{\lambda^r}{(\lambda L)^{i_0}}.\]
\item[(4)] For every step $1\leq s\leq i_0$ 
\[\|a^{(s+1)}-a^{(s)}\|_r\leq C\frac{\lambda^r}{(\lambda L)^{s}}\]
\end{itemize}
provided that, iteratively, the following inequalities hold for the bilinear error
\[\|R(a^{(s+1)},a^{(s+1)}-a^{(s)})\|_r+\|R(a^{(s)}-a^{(s+1)},a^{(s)})\|_r\leq \frac{C_r\lambda^r}{(\lambda \textcolor{black}{L})^{s+1}}.\]
\end{proposition}
\textsc{Remark:} the constant $C'>0$ might be large and depend on $i_0$, $C_*$, $k_1$ which are fixed. A proof of this an a generalization might be found in \cite{Mkroner}.

The following is a direct application of this proposition and Lemma \ref{conformal}.

\begin{lemma}[Final coefficients]\label{coefficients}
Under the hypothesis of Lemma \ref{conformal} and \ref{h}. For any step $s\geq s_0(\beta)$ there exist $A^{(s)}$ and constants $C>0$ such that the following holds:
\begin{itemize}
\item[(i)] $\|A^{(s)}_k\|_0\leq C$.
\item[(ii)] For any $1\leq r\leq \textcolor{black}{2}$
\[\|A^{(s)}_k \|_r\leq C\lambda_{q+1}^{r-1}\ell^{-1}.\]
\item[(iii)] Furthermore, for this choice of coefficients
\[\|h-M(A^{(s)})-L(A^{(s)})-R_5(A^{(s)})-R(A^{(s)})\|_r=O\left(\frac{\lambda_{q+1}^{r}}{(\lambda_{q+1}\ell)^{s}}\right).\]
\end{itemize}
\end{lemma}

\textsc{Remark:} the constants depend on $\epsilon$, $\theta$, $r_0$, the initial data $f_0$ and on the number of steps $s$ taken in the iteration but not on  $q$.

\textsc{Choice of parameters:} to prove (i) and (ii) we will apply K\"all\'en's iteration (i.e. Proposition \ref{kroner}) with the following choice of parameters: $r_1=2$, $i_0=s$ large (to be determined), $T=h$, $T_0=\operatorname{Id}$, $a^{(1)}=A^{(1)}$ from Lemma \ref{conformal}, $\lambda=\lambda_{q+1}$, $L=\ell$, $C_0>0$ satisfying $\sigma_1^{-1}\leq 3C_*$ large enough, the map $F$ will be defined as the inverse map corresponding to 
\[b(a,a)=M(a,a)+L(a)+R_5(a).\]
 Using the map $F$ given by Lemma \ref{decomposition} it can be written as
\[F(T)=F(T,\tau_k,\tau_{kk'})\]
noticing that, with this definition,
\[b(a,a)=\sum_{k=0}^{\frac{n(n+1)}{2}-1}a_k^2n_k\otimes n_k\textcolor{black}{+\sum_{k=0}^{\frac{n(n+1)}{2}-1}a_k\tau_k+\sum_{k,k'=0}^{\frac{n(n+1)}{2}-1}a_ka_{k'}\tau_{kk'}}\]
with $\tau_k$ and $\tau_{kk'}$ as defined in Lemma \ref{tau} (notice that this has been implicitly used in Lemma \ref{conformal} above). Taking into account the estimates from Lemma \ref{tau} together with Lemmata \ref{decomposition} and \ref{H4} we immediately obtain the required conditions on $F$. Notice that the constant $C_*$ is already fixed by $F$ and $\sigma_1$. We finally define the  error term $R$ by the formula
\[R(a)=R_2(a)+R_3(a)+R_4(a).\]
Before continuing any further let us introduce some back of the envelope calculations that shall be used in the proof. We will use these to check the initial step and the hypothesis of Proposition \ref{kroner} on the tensor and the error $R$. We will conclude proving part (iii) which is our main goal for the application.

\textsc{Back of the envelope calculations:} in this section we summarize some of the inequalities that hold in the appropiate range of parameters. We have already used that $\lambda_q\ell\leq 1$ which is a consequence of equation \ref{ell}. Furthermore, for any constant $C>0$ there is an appropiate choice in the double exponential ansatz such that
\begin{equation}\label{bound2}
\lambda_q\ell\leq 1\leq C\leq \ell\lambda_{q+1}.
\end{equation} 
Notice that the second condition implies $\lambda L$ can be chosen to be large enough which is an hypothesis on Proposition \ref{kroner} above. Indeed, we can ensure this by choosing $a$ large enough depending on $C_*$. Indeed, the condition phrased in this way is equivalent to $(\lambda_{q+1}\ell)^{-1}=o(1)$. Since we will need it at the end of the proof let us observe now a stronger statement, namely 
\begin{equation}\label{bound4}
(\lambda_{q+1}\ell)^{-s}=o(\delta_{q+2}\lambda_{q+1}^{-\epsilon}).
\end{equation}
For instance if $\epsilon=0$ and $s$ is large enough this essentially reduces to $2\alpha<b/\beta^*$. Indeed, unravelling $\ell$ it is equivalent to
\[\lambda_q^{s(1+\epsilon^*)}\delta_q^{\beta^*s}=o(\lambda_{q+1}^{s-\epsilon}\delta_{q+1}^{b+\beta^*s}).\]
Under the double exponential ansatz this follows provided
\[s(1+\epsilon^*)b^q-2\alpha\beta^*sb^{q-1}<(s-\epsilon)b^{q+1}-2\alpha(b+\beta^*s)b^q\]
which can be shown to imply
\[2\alpha<\frac{s(b-1)b}{b^2+\beta^*s(b-1)}+\epsilon\frac{b^2-s\beta^*b}{b^2+\beta^*s(b-1)}\]
It is now evident that for a fixed $b>1$ this restriction is irrelevant for small enough $\epsilon=\epsilon(b)>0$ when compared with the main one for a fixed but large enough $s\geq s_0(b,\beta)$   (cf. equation \ref{restriction2}).

It will be convenient to observe also that
\begin{equation}\label{bound0}
\lambda_q^{\beta-\epsilon}\ell^{-r+\beta}\leq C\frac{\lambda_{q+1}^r}{\lambda_{q+1}\ell}
\end{equation}
holds for any $r\geq 1$. This can be deduced from \ref{bound2} above. Another useful inequality will be
\begin{equation}\label{bound3}
\delta_q^{1/2}\lambda_q L\leq\delta_{q+1}^{1/2}\lambda L=\delta_{q+1}^{1/2}\lambda_{q+1}\ell\leq 1
\end{equation}
which follows provided
\[\alpha\left(1-2\beta^*(1-b^{-1})\right)\geq b-1-\epsilon^*\]
This follows for any $\epsilon>0$ provided it holds for $\epsilon=0$. The latter is true for any $\alpha>0$ provided $b=b(\alpha)$ is sufficiently close to one.

\textsc{Hypothesis on the tensor:} we can impose $\|T-T_0\|_0<(3C_*)^{-1}$ by choosing $\theta$ small enough (which can be done by letting $a$ be large enough in Lemma \ref{h}). Furthermore, the bounds for $T$, which involve $r\geq 1$, follow from Lemma \ref{h} and equation \ref{bound0}.

\textsc{The initial step $s=1$ of Proposition \ref{kroner}:} it is now evident that the initial step holds for (i), (ii) immediately from Lemma \ref{conformal} (parts (i) and (ii) respectively). \textcolor{black}{We are only left to check condition (iii) from  Proposition \ref{kroner} is satisfied. By the construction of $A^{(1)}$ in Lemma \ref{conformal} we have the identity $E_1=R(A^{(1)})$. Let us observe now the following straightforward inequalities hold:
\[\|R_2(A^{(1)})\|_r=\frac{O(1)}{\lambda_{q+1}^2}(\delta_q^{1/2}\lambda_q)^2\lambda_{q+1}^r,\]
\[\|R_3(A^{(1)})\|_r=\frac{O(1)}{\lambda_{q+1}^2}\left(\frac{\lambda_{q+1}}{\lambda_{q+1}\ell}\right)^2\lambda_{q+1}^r\]
and\[\|R_4(A^{(1)})\|_r=\frac{O(1)}{\lambda_{q+1}^2}\frac{\lambda_{q+1}}{\lambda_{q+1}\ell}\delta_q^{1/2}\lambda_q\lambda_{q+1}^r\]
as can be checked by the reader by inspecting the definition of $R_2$, $R_3$ and $R_4$ above. The bound $\|R\|_r=O(\lambda^r(\lambda L)^{-1})$ follows because it holds true for each of the above error terms individually. This can be easily checked as a consequence of equation \ref{star}.}

\textsc{The bilinear term of type $R_2$:} using the definition and the \textcolor{black}{available bounds on the derivatives of the normals $\partial_i \nu$ (cf. Lemma \ref{basis})} we can estimate
\[\|R_2(a,b)\|_r\leq C_r\frac{\delta_q\lambda_q^2}{\lambda_{q+1}^2}\sum_{j_1+j_2=0}^r\|a\|_{j_1}\|b\|_{j_2}\lambda_{q+1}^{r-j_1-j_2}\]
which, in the case of evaluating at \textcolor{black}{$\{a,b\}=\{A^{(s)}, A^{(s)}-A^{(s+1)}\}$ and $\{a,b\}=\{A^{(s+1)}, A^{(s)}-A^{(s+1)}\}$} using that the difference satisfy Proposition \ref{kroner} (4) we get the bound $O(\lambda^{r}(\lambda L)^{-s-1})$ immediately using equation \ref{bound3}. 

\textsc{The bilinear term of type $R_3$:}  using the definition and the bounds for the normals we can bound
\[\|R_3(a,b)\|_r\leq C_r\sum_{j_1+j_2=0}^r\|a\|_{j_1+1}\|b\|_{j_2+1}\lambda_{q+1}^{r-j_1-j_2-2}\]
which, in the case of evaluating with \textcolor{black}{$\{a,b\}=\{A^{(s)}, A^{(s)}-A^{(s+1)}\}$ and $\{a,b\}=\{A^{(s+1)}, A^{(s)}-A^{(s+1)}\}$} using that the difference satisfy Proposition \ref{kroner} (4) and (ii) we get the bound $O(\ell^{-2}\lambda^{r-2}(\lambda L)^{-s-1})$. 

\textsc{The bilinear term of type $R_4$:} in this case we need to be slightly more careful than in the previous. Again, using the definition and the bounds for the normals we can bound
\[\|R_4(a,b)\|_r\leq C_r\delta_{q}^{1/2}\lambda_q\sum_{j_1+j_2=0}^r\|a\|_{j_1}\|b\|_{j_2+1}\lambda_{q+1}^{r-j_1-j_2-2}\]
which, in the case of evaluating with \textcolor{black}{$\{a,b\}=\{A^{(s)}, A^{(s)}-A^{(s+1)}\}$ and $\{a,b\}=\{A^{(s+1)}, A^{(s)}-A^{(s+1)}\}$} using that the difference satisfy Proposition \ref{kroner} (4) and equation \ref{bound3} we get the bound $O(\lambda^{r}(\lambda L)^{-s-1})$ as before. 

This ensures the induction that allows to use K\"all\'en's iteration  is closed but to conclude the proof we still need to check \textcolor{black}{(iii)}.

\textsc{End of proof, part} (iii): by definition 
\[b(A^{(s)})=M(A^{(s)})+L(A^{(s)})+R_5(A^{(s)})\] 
which gives
\[h-b(A^{(s)})-R(A^{(s)})=h-M(A^{(s)})-L(A^{(s)})-R_5(A^{(s)})-R(A^{(s)})=E_s\]
which corresponds to the error in estimate (3) in Proposition \ref{kroner}.

The following result merges all the previous to provide the perturbation we were seeking for.

\begin{lemma}[Perturbation]\label{perturbation}
Under the hypothesis of Proposition \ref{induction}. There exist a perturbation $w_{q+1}$ such that $\|w_{q+1}\|_r\leq C\delta_{q+1}^{1/2}\lambda_{q+1}^{r-1}$ for $r\in[0,2]$,
\[\|g(x,y)-g^{q+1}(x,y)-\delta_{q+2}\operatorname{Id}\|_0\leq\delta_{q+2}\lambda_{q+1}^{-\epsilon}\]
and
\[\|g(x,y)-g^{q+1}(x,y)-\delta_{q+2}\operatorname{Id}\|_{\beta}\leq\delta_{q+2}\lambda_{q+1}^{\beta-\epsilon}\]
\end{lemma}

\textsc{Proof:} we will produce now a perturbation as indicated in Section \ref{theperturbation} with a choice of coefficients $A$ given by $A^{(s)}$ from Lemma \ref{coefficients}. By definition and the bounds in Lemma \ref{coefficients} one observes 
\[\|w_{q+1}\|_0\leq C\delta_{q+1}^{1/2}\lambda_{q+1}^{-1}\] 
 On the other hand by interpolation (cf. Section \ref{sectionholder}) it is enough to prove that
\[\|w_{q+1}\|_1\leq C\delta_{q+1}^{1/2}\textrm{ and }\|w_{q+1}\|_2\leq C\delta_{q+1}^{1/2}\lambda_{q+1}\]
to gain the required control of the perturbation in the full range $r\in[0,2]$. To prove these we might observe first that
\[\|w_{q+1}\|_1=\sum_{k=0}^{\frac{n(n+1)}{2}-1}\frac{\delta_{q+1}^{1/2}}{\lambda_{q+1}}O\left(\|A\|_0\|\nu_k\|_1+\|A\|_1\|\nu_k\|_0+\lambda_{q+1}\|A\|_0\|\nu_0\|_0\right)\]
Using this, together with the previous observation and Lemmata \ref{basis} and \ref{conformal} we conclude
\[\|w_{q+1}\|_1=\frac{\delta_{q+1}^{1/2}}{\lambda_{q+1}}O(\delta_q^{1/2}\lambda_q+\ell^{-1}+\lambda_{q+1})\leq C\delta_{q+1}^{1/2}\]
which follows taking into account inequality \ref{bound2}. Analogously, expanding two derivatives of $w_{q+1}$, from its definition above, the terms where no derivatives are attached to $\nu_k$ are bounded by
\[\sum_{k=0}^{\frac{n(n+1)}{2}-1}\frac{\delta^{1/2}_{q+1}}{\lambda_{q+1}}\left(\|A\|_2+\lambda_{q+1}\|A\|_1+(\lambda_{q+1})^2\|A\|_0\right).\]
On the other hand if any derivative hits $\nu_i$ one gets the bound
\[\sum_{k=0}^{\frac{n(n+1)}{2}-1}\frac{\delta_{q+1}^{1/2}}{\lambda_{q+1}}O(\|\nu_k\|_2\|A\|_0+\|\nu_k\|_1(\|A\|_1+\lambda_{q+1}\|A\|_0)).\]
We deal with them separately. The first can be bounded using Lemma \ref{coefficients} by
\[\frac{\delta_{q+1}^{1/2}}{\lambda_{q+1}}O(\lambda_{q+1}\ell^{-1}+\lambda_{q+1}\ell^{-1}+(\lambda_{q+1})^2)\]
where the constant depends on the bounds from Lemma \ref{coefficients} and, in particular, is independent of $q$. Again, using equation \ref{bound2} this bound simplifies to $O(\delta_{q+1}^{1/2}\lambda_{q+1})$. Using the same strategy, Lemmata \ref{basis} and \ref{coefficients}, the second can be bounded as follows
\[\frac{\delta_{q+1}^{1/2}}{\lambda_{q+1}}O\left(\delta_q^{1/2}\lambda_q\ell^{-1}+\delta_q^{1/2}\lambda_q(\ell^{-1}+\lambda_{q+1})\right)=O(\delta_{q+1}^{1/2}\lambda_{q+1})\]
which proves our claim for $r=2$. The rest follows by interpolation.

Notice that, since we defined $f_{q+1}=f^{\ell}_q+w_{q+1}$, the above implies 
\[\|f_{q+1}\|_r\leq\|f_q\|_r+\|w_{q+1}\|_r\leq \kappa_q^r+C\delta_{q+1}^{1/2}\lambda^{r-1}_{q+1}.\]

Finally observe that with our notational convention (cf. Proposition \ref{kroner} and Lemma \ref{coefficients})
\[\begin{aligned}
g_{ij}-g_{ij}^{q+1}&=g_{ij}-(g^q_{\ell})_{ij}-2\partial_if^{\ell}_q\cdot\partial_jw_{q+1}-\partial_iw_{q+1}\cdot\partial_jw_{q+1}\\
&=\delta_{q+1}(E_{s}+R_6).
\end{aligned}\]
Taking into account the cancellation $R_6=0$ and the bounds given by Lemma \ref{coefficients} (iv) for $E_s$ one obtains
\[\|g(x,y)-g^{q+1}(x,y)-\delta_{q+2}\textrm{Id}\|_r\leq\delta_{q+1}\frac{C(s_0)\lambda_{q+1}^r}{(\lambda_{q+1}\ell)^{s_0}}\]
with $r=0,1$. The error is of order $o(\delta_{q+2}\lambda_{q+1}^{r-\epsilon})$ due to equation \ref{bound4} above. Finally, interpolating produces the required $C^{\beta}$ norm and finishes the proof of Proposition \ref{induction}.

\section{The initial step and proof of Theorem \ref{first}}\label{prooffirst}

In this section we end the proof of the main result by providing an initial step for the induction. Without loss of generality we can assume the short embedding $h$ is such that $g-h^{\sharp}e-\eta \textrm{Id}> 0$ pointwisely in the quadratic forms sense for some $\eta$. By taking $\eta$ even smaller if necessary we can also suppose  that $\eta=\delta_1=a^{-2\alpha }$ for some $a$ large enough. An application of the refinement of the Nash-Kuiper theorem (Theorem \ref{nashkuiper}) provides the existence of a $C^{1,\alpha_0}$ embedding $h_0$ such that such that $g-h_0^{\sharp}e-\eta\textrm{Id}=0$ and $\|h-h_0\|_{0}\leq\epsilon$. To fulfill the hypothesis of Proposition \ref{induction} we need to control higher norms. This suggest to use a small mollification (at some small scale $\ell$) of $h_0$, say  $f_0$, so that $\|f_0\|_2$ is bounded by some (maybe large) constant and arguing as in the proof of Lemma \ref{mollification} (cf. equation \ref{moleq} where the middle term vanishes by our construction) one gets
\[\left\|\frac{g-f_0^{\sharp}e}{\eta}-\textrm{Id}\right\|_{\beta}\leq C\eta^{-1}(\ell^2\|g_{ij}\|_{2, \beta}+\ell^{2\alpha_0-\beta}\|h_0\|_{1,\alpha_0}^2)\]
One can take now $\ell$ very small so that this is less than $a^{-1}=\lambda_0^{-1}<\sigma_1$. Notice that this requires $\beta<2\alpha_0$ which together with equation \ref{restriction2} provides the limitation $\beta<1$. 

The induction step is precisely Proposition \ref{induction} which produces a sequence of functions $f_q$. The limit would be well defined in $C^{1,1-\epsilon'}$ if we can bound the sequence uniformly in a space where it is compactly embedded (cf. Lemma \ref{compact}). Indeed, as a consequence of the bounds in the Proposition \ref{induction}  for $2\alpha<2-\beta$ one gets 
\[\|f_q\|_{1,\alpha'}\leq\|f_0\|_{C^{1,\alpha'}}+C\sum_{q=1}^{\infty}\delta_q^{1/2}\lambda_q^{\alpha'}\]
which is bounded thanks to the double exponential ansatz for any $\alpha'<\alpha/b$. This finally shows the existence of a subsequence which converges to  an isometric embedding $f$ in $C^{1,1-\epsilon'}$ for any $\epsilon'>0$ provided we choose the parameters $\beta$ small enough, $b$ close enough to one and $a$ large enough. Furthermore, by taking $\ell=\ell(\epsilon)$ small enough and $a$ large enough we can make
\[\|h-f\|_{0}\leq \|h-h_0\|_0+\|h_0-f_0\|_0+\sum_{q=1}^{\infty}\delta_q^{1/2}\lambda_q^{-1}<3\epsilon.\]
Indeed, the first is bounded by $\epsilon$ by Theorem \ref{nashkuiper} and our construction, the second is bounded by $C\ell\|h_0\|_1$ by Lemma \ref{holder} while the last one is as small as we please choosing $a$ large enough.

\textsc{Remark:} the reader might want to compare this with De Lellis, Inauen and Sz\'ekelyhidi who construct their initial $f_0$ using two approximations in their paper \cite{DIS}. Observe also that the proof provides even more that stated in Theorem \ref{first}, namely $g_q\rightarrow g$ in $C^{\epsilon}$.

\section{The perturbation using strains and proof of Theorem \ref{first}}\label{theperturbationone}

The perturbation $w_{q+1}$ can be adapted to deal with the codimension one case. In this section we sketch how to adapt the method. In this case, we follow \textcolor{black}{closely} the original construction by Kuiper, using a \textcolor{black}{new} variant of strains instead of spirals (cf. \cite{K}, sections 6 and 7; an alternative to strains are corrugations, see \cite{CLS}). 

We define $f^{(0)}=f_q^{\ell}$ and $f^{(k)}=f^{(k-1)}+w_{q+1}^{(k-1)}$ for $k\geq 1$ as before. We will define $\nu^1_k(x)$ as the  directional derivative $\nabla \textcolor{black}{f_q^{\ell}}(x)\cdot v_k$ where 
\[v_k
^r=(g_{\ell}^q)^{ir}(n_k)_i.\]
We are using the notation $g^{ir}$ to denote $(g^{-1})_{ir}$ as usual. On the other hand the $\nu^2_k$, $k=0,\ldots,\frac{n(n+1)}{2}-1$, will still correspond to an orthonormal basis of vectors orthogonal to the tangent plane of $f_q^{\ell}$ at $x$ (cf. Lemma \ref{basis} with $d=\frac{n(n+1)}{2}$ for a construction). Notice that, by definition, every $\nu^1$ is a linear combination of tangent vectors. In particular, they are orthogonal to the $\nu^2$. Let us digress on their properties a little bit more and prove the following
\begin{lemma}\label{nu}
The normals $\nu^1$ satisfy the bounds $\|\nu^1\|_0\leq C$ and
\[\|\nu^1\|_r\leq C_r\ell^{-r}\]
for $r\geq 1$.
\end{lemma}
\textsc{Remark:} notice that $\nu^1$ depends implicitly on $q$ while the constant depends on the fixed metric $g$ but not on the step $q$. 

\textsc{Proof:} their definition requires to invert the tensor $g_{\ell}^q$ which is close to the fixed metric $g$. Using Lemma \ref{mollification} it is evident that it satisfies the bounds
\[\|g_{\ell}^q\|_0\leq C(g)\]
and
\[\|g_{\ell}^q\|_r\leq C(g)+\delta_{q+1}\lambda_q^{-\epsilon}(\lambda_q\ell)^{\beta}\ell^{-r}\]
for $r\geq 1$. The inverse of $g_q^{\ell}$ exists since it exist for $g$ and  also for any other metric in the  neighbourhood of $g$. That $g_q^{\ell}$ is close enough to $g$ can be fulfilled if $\delta_1$ is small enough for a fixed $g$, using Lemma \ref{mollification}. Actually, by continuity,  reducing the neighbourhood if necesary, we can impose their determinants to be comparable. It is now easy, observing that the inverse can be written algebraically using the functions defining the mollified tensor in a canonical basis, to show that
\[\|(g_{\ell}^q)^{-1}\|_0\leq C(g)\]
and
\[\|(g_{\ell}^q)^{-1}\|_r\leq C(g)(1+\delta_{q+1}\lambda_q^{-\epsilon}(\lambda_q\ell)^{\beta}\ell^{-r}).\]
Using these estimates, Leibniz's rule, the estimates for $f_q^{\ell}$ from Lemma \ref{mollification} and the inequality $\delta_{q+1}\lambda_q^{-\epsilon}(\lambda_q\ell)^{\beta}=O(1)$ it follows that
\[\|\nu^1\|_r\leq C_r\left(\delta_q^{1/2}\lambda_q+\ell^{-1}\right)\ell^{-r+1}\]
which implies our claim taking into account equation \ref{bound3}.

The reader might want to check that these estimates are consistent with the proof of Lemma \ref{perturbation} above. On the other hand, the $\nu^2$ satisfy identical estimates as our previous $\nu$. The normals of type $\nu^1_k$ are not orthogonal to the tangent plane but satisfy the following fundamental identity
\begin{equation}\label{n}
\partial_if_q^{\ell}\cdot \nu^1_k=(n_k)_i.
\end{equation}
Its proof is almost immediate from the definition. Indeed, for any $i=1,\ldots,n$ we have
\[\partial_if\cdot \nu^1_k=\sum_{r,s=1}^n\partial_if_s\partial_rf_s v^r_k=\sum_{r=1}^ng_{ir}v^r_k=\sum_{r,t=1}^ng_{ir}g^{tr}(n_k)_t=(n_k)_i\]
where $f=f_q^{\ell}$ and $g=g^q_{\ell}$. This time, the perturbations will be based on \textit{strains} having the form
\[\textcolor{black}{w^{(k)}_{q+1}(x)=\frac{\delta_{q+1}A_k(x)^2}{4\lambda_{q+1}}\sin(2\lambda_{q+1}x\cdot n_{k})\nu^1_{k}(x)+\delta^{1/2}_{q+1}\textcolor{black}{\sqrt{2}}\frac{A_k(x)}{\lambda_{q+1}}\cos\left(\lambda_{q+1}x\cdot n_{k}\right)\nu^2_{k}(x)}\]
where the frequencies are given by $\lambda_{q+1}=\lambda_q^{b}$. (This definition might be compared with equation 6.7 in \cite{K}.)  The choice of coefficients $A_k$ will be specified later \textcolor{black}{as in the proof of Theorem \ref{kroner} above}. This time we get
\protect{\small\[\begin{aligned}
g^{q+1}&=g_{\ell}^q+\sum_{k=0}^{\frac{n(n+1)}{2}-1}\partial_i\textcolor{black}{f_q^{\ell}}\cdot\partial_j\left(\delta_{q+1}\frac{A_k(x)^2}{4\lambda_{q+1}}\sin(2\lambda_{q+1}x\cdot n_{k})\right)\nu^1_{k}(x)\\
&\qquad\qquad+\partial_i\textcolor{black}{f_q^{\ell}}\cdot\partial_j\left(\delta_{q+1}^{1/2}\textcolor{black}{\sqrt{2}}\frac{A_k(x)}{\lambda_{q+1}}\cos\left(\lambda_{q+1}x\cdot n_{k}\right)\right)\nu^2_{k}(x)\\
&+\sum_{k=0}^{\frac{n(n+1)}{2}-1}\partial_i\textcolor{black}{f_q^{\ell}}\cdot\left(\partial_j\nu^1_{k}(x)\frac{\delta_{q+1}A_k^2(x)}{4\lambda_{q+1}}\sin(2\lambda_{q+1}x\cdot n_{k})\right.\\
&\left.\qquad\qquad\qquad+\partial_j\nu^2_{k}(x)\delta_{q+1}^{1/2}\textcolor{black}{\sqrt{2}}\frac{A_k(x)}{\lambda_{q+1}}\cos\left(\lambda_{q+1}x\cdot n_{k}\right)\right)\\
&+\sum_{k\textcolor{black}{,k'}=0}^{\frac{n(n+1)}{2}-1}\partial_iw^{(k)}_{q+1}\cdot\partial_jw^{(\textcolor{black}{k'})}_{q+1}
\end{aligned}\]}
\textcolor{black}{Again, we are abusing the notation since two more terms corresponding to the second and third summands with $i$ and $j$ reversed should appear.} It will be convenient to introduce the following notation for a relevant \textcolor{black}{analogue to our previous {\em linear} term}
\[\begin{aligned}
L(A)&=\sum_{k=0}^{\frac{n(n+1)}{2}-1}\partial_if_q^{\ell}\cdot\left(\partial_j\nu^1_{k}(x)\frac{A_k^2(x)}{4\lambda_{q+1}}\sin(2\lambda_{q+1}x\cdot n_{k})\right.\\
&\left.+\partial_j\nu^2_{k}(x)\textcolor{black}{\sqrt{2}}\frac{A_k(x)}{\lambda_{q+1}\delta_{q+1}^{1/2}}\cos\left(\lambda_{q+1}x\cdot n_{k}\right)\right).
\end{aligned}\]
Let us introduce the following linear operator acting on tensors
\[(\operatorname{sym}A)_{ij}=\frac{1}{2}(A_{ij}+A_{ji}),\]
that leaves invariant symmetric tensors. We will use it in the sequel to correct our previous abuse of notation. Then we get that
\[\begin{aligned}
g^{q+1}&=g^q_{\ell}+\operatorname{sym}\sum_{k=0}^{\frac{n(n+1)}{2}-1}\partial_i\textcolor{black}{f_q^{\ell}}\cdot\partial_j\left(\delta_{q+1}\frac{A_k(x)^2}{2\lambda_{q+1}}\sin(2\lambda_{q+1}x\cdot n_{k})\right)\nu^1_{k}(x)\\
&\qquad+2\delta_{q+1}\operatorname{sym}L+\sum_{k\textcolor{black}{,k'}=0}^{\frac{n(n+1)}{2}-1}\partial_iw^{(k)}_{q+1}\cdot\partial_jw^{(\textcolor{black}{k'})}_{q+1}
\end{aligned}\]
because the normals $\nu^2$ are orthogonal to the tangent. We should remark now that, in this case, part of the metric error is approximated using the linear part too. This justifies the introduction of the following main term
\begin{equation}\label{mainterm}
\begin{aligned}
M(a)_{ij}&=\operatorname{sym}\sum_{k=0}^{\frac{n(n+1)}{2}-1}\partial_i\textcolor{black}{f_q^{\ell}}\cdot a_k(x)^2\cos(2\lambda_{q+1}x\cdot n_{k})\nu^1_{k}(x)(n_k)_j\\
&\textcolor{black}{\qquad+2\sum_{k,k'=0}^{\frac{n(n+1)}{2}-1}a_k(x)\sin\left(\lambda_{q+1}x\cdot n_{k}\right)\nu^2_{k}(x) \cdot}\\
&\qquad\qquad\cdot a_{k'}(x) \sin\left(\lambda_{q+1}x\cdot n_{k'}\right)\nu^2_{k'}(x)(n_k)_i (n_{k'})_j\\
&=\sum_{k=0}^{\frac{n(n+1)}{2}-1}a_k(x)^2(n_k\otimes n_{\textcolor{black}{k}})_{ij}
\end{aligned}
\end{equation}
which can be justified using the orthogonality of the normals, equation \ref{n} and the elementary identity
\[\cos(2z)=1-2\sin^2(z).\]
Let us introduce some notation analogous to our previous analysis for the relevant tensor valued error terms. The first that we list, 
\textcolor{black}{\protect{\small \[R_0(a)_{ij}=\sum_{k,k'=0}^{\frac{n(n+1)}{2}-1}\delta_{q+1}^{1/2}\frac{a_k(x)^2}{2}\sqrt{2}a_{k'}(x) \cos(2\lambda_{q+1}x\cdot n_k)\sin(\lambda_{q+1}x\cdot n_{k'})\nu_k^1(x)\cdot\nu_{k'}^2(x)(n_k\otimes n_k)_{ij}=0,\]
vanishes by orthogonality. It is, nevertheless, not bilinear. This will also happen with the rest of terms. The next three  are symmetric tensors:
\[R_1(a)_{ij}=\sum_{k,k'=0}^{\frac{n(n+1)}{2}-1}\delta_{q+1}\frac{a_k(x)^2}{2}\frac{a_{k'}(x)^2}{2}\cos(2\lambda_{q+1}x\cdot n_{k})\nu^1_{k}(x) \cdot\cos(2\lambda_{q+1}x\cdot n_{k'})\nu^1_{k'}(x)(n_k\otimes n_k)_{ij},\]
\[\begin{aligned}
R_2(a)_{ij}=\sum_{k,k'=0}^{\frac{n(n+1)}{2}-1}&\frac{1}{\lambda_{q+1}^2}\left(\frac{\delta_{q+1}^{1/2}}{4}a_k(x)^2\sin(2\lambda_{q+1}x\cdot n_{k})\partial_i\nu^1_{k}(x)+  \textcolor{black}{\sqrt{2}}a_{k}(x)\cos(\lambda_{q+1}x\cdot n_{k})\partial_i\nu^2_{k}(x)\right) \cdot\\
 & \cdot\left(\frac{\delta_{q+1}^{1/2}}{4}a_{k'}(x)^2\sin(2\lambda_{q+1}x\cdot n_{k'})\partial_j\nu^1_{k'}(x)+\textcolor{black}{\sqrt{2}}a_{k'}(x)\cos(\lambda_{q+1}x\cdot n_{k'})\partial_j\nu^2_{k'}(x)\right),
\end{aligned}\]
and
\[\begin{aligned}
R_3(a)_{ij}=\sum_{k,k'=0}^{\frac{n(n+1)}{2}-1}&\frac{1}{\lambda_{q+1}^2}\left(\frac{\delta_{q+1}^{1/2}}{4}\partial_ia_k(x)^2\sin(2\lambda_{q+1}x\cdot n_{k})\nu^1_{k}(x)+\textcolor{black}{\sqrt{2}}\partial_ia_{k}(x)\cos(\lambda_{q+1}x\cdot n_{k})\nu^2_{k}(x) \right) \cdot\\
 & \cdot\left(\frac{\delta_{q+1}^{1/2}}{4}\partial_ja_{k'}(x)^2\sin(2\lambda_{q+1}x\cdot n_{k'})\nu^1_{k'}(x)+\textcolor{black}{\sqrt{2}}\partial_ja_{k'}(x)\cos(\lambda_{q+1}x\cdot n_{k'})\nu^2_{k'}(x)\right).
\end{aligned}\]
The rest of them are not symmetric, namely:
\[\begin{aligned}
R_4(a)_{ij}=\sum_{k,k'=0}^{\frac{n(n+1)}{2}-1}&\frac{1}{\lambda_{q+1}^2}\left(\frac{\delta_{q+1}^{1/2}}{4}a_k(x)^2\sin(2\lambda_{q+1}x\cdot n_{k})\partial_i\nu^1_{k}(x)+  \textcolor{black}{\sqrt{2}}a_k(x)\cos(\lambda_{q+1}x\cdot n_{k})\partial_i\nu^2_{k}(x)\right) \cdot\\
 & \cdot\left( \frac{\delta_{q+1}^{1/2}}{4}\partial_ja_{k'}(x)^2\sin(2\lambda_{q+1}x\cdot n_{k'})\nu^1_{k'}(x)+\textcolor{black}{\sqrt{2}}\partial_ja_{k'}(x)\cos(\lambda_{q+1}x\cdot n_{k'})\nu^2_{k'}(x)\right),
\end{aligned}\]
\[\begin{aligned}
R_5(a)_{ij}=\sum_{k,k'=0}^{\frac{n(n+1)}{2}-1}&\frac{1}{\lambda_{q+1}}\left(\frac{\delta_{q+1}^{1/2}}{4}a_k(x)^2\sin(2\lambda_{q+1}x\cdot n_{k})\partial_i\nu^1_{k}(x)+  \textcolor{black}{\sqrt{2}}a_k(x)\cos(\lambda_{q+1}x\cdot n_{k})\partial_i\nu^2_{k}(x)\right) \cdot\\
 & \cdot\left(\frac{\delta_{q+1}^{1/2}}{2}a_{k'}(x)^2\cos(2\lambda_{q+1}x\cdot n_{k'})\nu^1_{k'}(x)-\textcolor{black}{\sqrt{2}}a_{k'}(x)\sin(\lambda_{q+1}x\cdot n_{k'})\nu^2_{k'}(x)\right)(n_{k'})_j,
\end{aligned}\]
\[\begin{aligned}
R_6(a)_{ij}&=\sum_{\textcolor{black}{k,k'=0}}^{\frac{n(n+1)}{2}-1}\frac{1}{\lambda_{q+1}}\left(\frac{\delta_{q+1}^{1/2}}{4}\partial_ia_k(x)^2\sin(2\lambda_{q+1}x\cdot n_{k})\nu^1_{k}(x)+ \textcolor{black}{\sqrt{2}}\partial_ia_k(x)\cos(\lambda_{q+1}x\cdot n_{k})\nu^2_{k}(x) \right) \cdot\\
 & \qquad\qquad\cdot\left(\frac{\delta_{q+1}^{1/2}}{2}a_{k'}(x)^2\cos(2\lambda_{q+1}x\cdot n_{k'})\nu^1_{k'}(x)-\textcolor{black}{\sqrt{2}}a_{k'}(x)\sin(\lambda_{q+1}x\cdot n_{k'})\nu^2_{k'}(x)\right)(n_{k'})_j
\end{aligned}\]
and, finally,
\[R_7(a)_{ij}=\sum_{k=0}^{\frac{n(n+1)}{2}-1}\partial_if_q^{\ell}\cdot\frac{a_k(x)\partial_j a_k(x)}{2\lambda_{q+1}}\sin(2\lambda_{q+1}x\cdot n_{k})\nu^1_{k}(x).\]}}

We observe that there are new error terms, namely $R_1$ and $R_7$. Furthermore, if we unfold the products it becomes evident that this error terms contain bilinear, trilinear and cuatrilinear expressions in $a$ (cf. $R_1$ above). Let us denote by $\tilde{R}_i$ the trilinear and cuatrilinear parts of $R_i$. As before one can readily check that
\[\begin{aligned}
g^{q+1}&=g_{\ell}^q+\delta_{q+1}\operatorname{sym}(2L+M+R_1+R_2+R_3+R_4+R_5+R_6+2R_7)\\
&=g_{\ell}^q+\delta_{q+1}(\operatorname{sym}(2L+M+R_1+R_2+R_3+R_4+R_5)+\operatorname{sym}(R_6+2R_7)).
\end{aligned}\]
This formula can be used to observe that although $L$, $R_4$, $R_5$ and $R_6$ and $R_7$ are not symmetric by definition, but  jointly they are. The specific grouping of terms will be seen to be relevant later.

\textcolor{black}{We might then try to apply Proposition \ref{kroner} as before to get rid of the extra error terms. Before proceeding any further let us digress and introduce the following substitute for Lemma \ref{tau} above
\begin{lemma}\label{tau2}
\textcolor{black}{The symmetric tensors
\[\begin{aligned}
(\tau_k)_{ij}&=2\operatorname{sym}\partial_if_q^{\ell}\cdot\partial_j\nu^2_{k}(x)\frac{\textcolor{black}{\sqrt{2}}}{\lambda_{q+1}\delta_{q+1}^{1/2}}\cos\left(\lambda_{q+1}x\cdot n_{k}\right)\\
\end{aligned}\]
and 
\[\begin{aligned}
(\tau_{kk'})_{ij}&=-\frac{1}{\lambda_{q+1}}\textcolor{black}{\sqrt{2}}\operatorname{sym}\cos(\lambda_{q+1}x\cdot n_{k})\partial_i\nu^2_{k}(x)\cdot\textcolor{black}{\sqrt{2}}\sin(\lambda_{q+1}x\cdot n_{k'})\nu^2_{k'}(x))(n_{k'})_j\\
 &\qquad+\delta_{kk'}\frac{1}{2\lambda_{q+1}}\operatorname{sym}\partial_if_q^{\ell}\cdot\partial_j\nu^1_{k}(x)\sin(2\lambda_{q+1}x\cdot n_{k}) 
 \end{aligned}\]
 satisfy the estimates
 \[\|\tau_k\|_r+\|\tau_{kk'}\|_r\leq\frac{\lambda_{q+1}^r}{\lambda_{q+1}\ell}.\]}
\end{lemma}
Notice that this tensors are intimately related the bilinear parts of $L$ and $R_5$. The proof of this result is identical to the one of Lemma \ref{tau} employing Lemma \ref{nu} to deal with the term involving $\nu^1$. We leave further details to the reader.}

Our goal now would be to provide a substitute for Lemma \ref{coefficients}. The bilinear parts of the errors of types $R_2$, $R_3$, $R_4$ and $R_5$ are virtually the same. Indeed, notice that they depend only on the bounds on $\nu^2$, which are the same as in the proof of Theorem \ref{kroner}. Nevertheless there are some new obstacles too. As already mentioned, there are trilinear and cuatrilinear error terms that we will need to handle separately, the error $R_6$ does not vanish nor is it suitable to K\"all\'en's iteration and there is a similar error term $R_7$ which is also not suitable to the iteration. 

Before continuing any further let us observe that $R_1$ can be kept out of the iteration. Indeed, notice that it can be bounded in $C^r$ by $O(\delta_{q+1}\lambda_{q+1}^r)$, it is then easy to observe that the corresponding metric error term $\delta_{q+1}R_1$ will close the induction. In general,
\begin{lemma}
The following estimate holds
\[\|\tilde{R}_i\|_r=O(\delta_{q+1}^{1/2}\lambda_{q+1}^r)\]
for any $i=1,2,3,4,5,6$ and $r\geq 0$.
\end{lemma} 
We leave checking this to the reader's discretion. This estimates  allow to also leave the rest of trilinear and cuatrilinear parts out of K\"all\'en's iteration.

Our analysis reduces then to handle the errors of type $R_6$ and $R_7$ appropiately. A rough estimate for $R_6$ or $R_7$ would lead to the restriction $\alpha<1/3$ which is certainly weaker than the one claimed in the statement of Theorem \ref{first}.  The same issue happens with $R_7$. It is nevertheless surprising to observe that taken together they satisfy
\[\|(R_6)_{ij}+2(R_7)_{ji}\|_r=O\left(\frac{\delta_{q+1}\lambda_{q+1}^r}{\lambda_{q+1} \ell}\right).\]
This estimate follows as a consequence of an algebraic identity that cancels $R_7$ with the bilinear part of $R_6$. Indeed, one can show such a cancellation using equation \ref{n} and the elementary $\sin(2z)=2\sin(z)\cos(z)$ to prove
\textcolor{black}{\[\begin{aligned}
2(R_7)_{ij}=&\frac{1}{\lambda_{q+1}}\sum_{k=0}^{\frac{n(n+1)}{2}-1}a_k(x)\partial_ja_k(x)2\sin(\lambda_{q+1}x\cdot n_k)\cos(\lambda_{q+1}x\cdot n_k)(n_k)_i\\
&=(\tilde{R}_6-R_6)_{ji}.
\end{aligned}\]}
This allows us to deal with the error terms as in the proof of  Lemma \ref{perturbation}: either separately or within K\"all\'en's iteration (using a substitute for Lemma \ref{coefficients}). We would still need to check that the estimates in Lemma \ref{perturbation} for the perturbation $w_{q+1}$ hold. This can be done using Lemma \ref{nu} to take into account the new bounds for the normals $\nu^1$. We leave further details to the reader. 

 An extension of this following the setting of {\em adapted short immersions} introduced by Cao and Sz\'ekelyhidi, provides the global version of this result (cf. \cite{Ka, CS3}). 

\section{Appendix: H\"older spaces and mollification}\label{sectionholder}

For any integer $k\geq 0$ and $\alpha\in (0,1)$ the H\"older space $C^{k,\alpha}(B)$ is defined as the set of functions whose norm
\[\|f\|_{k,\alpha}=\|f\|_{k+\alpha}=\sup_{x\in B}\sum_{|a|\leq k}|\partial^af(x)|+\sup_{x\neq y\in B}\sum_{|a|=k}\frac{|\partial^af(x)-\partial^af(y)|}{|x-y|^{\alpha}}\]
is finite. Tensors are measured with respect to the usual H\"older norms in local coordinates componentwise. 

\begin{lemma}\label{compact}
Let  $n\geq 0$ and $\alpha<\beta$ then $Id:C^{n,\beta}\rightarrow C^{n,\alpha}$ is compact.
\end{lemma}

\begin{lemma}[Interpolation]
Given $\lambda\in (0,1)$ such that $\alpha=\lambda\alpha_1+(1-\lambda)\alpha_2$ the following inequality holds
\[\|f\|_{\alpha}\leq \|f\|_{\alpha_1}^{\lambda}\|f\|_{\alpha_2}^{1-\lambda}.\]
\end{lemma}

Let us introduce how we intend to mollify functions. Fix a positive smooth compactly supported symmetric $\varphi$ such that $\int\varphi=1$. We will need to mollify at different scales which justifies the introduction of $\varphi_{\ell}(x)=\ell^{-n}\varphi(x\ell^{-1})$. The new $\varphi_{\ell}$ are also positive smooth compactly supported symmetric function with unit mass. Given a function $f$ we define its mollification at scale $\ell$, usually denoted $f_{\ell}$, as $f*\varphi_{\ell}$. Of special importance for us will be the following

\begin{lemma}[cf. \cite{CLS, DIS}]\label{holder}
Let $\varphi\in C^{\infty}_c(\mathbb{R}^n)$ as explained before. Then for any $r,s\in\mathbb{N}$ and $\alpha\in (0,1]$ we have
\begin{itemize}
\item[(i)] $[f*\varphi_{\ell}]_{r+s}\leq C\ell^{-s}\|f\|_r$,
\item[(ii)] $[f-f*\varphi_{\ell}]_r\leq C\ell^2\|f\|_{r+2}$,
\item[(iii)] $[(fg)*\varphi_{\ell}-(f*\varphi_{\ell})(g*\varphi_{\ell})]_r\leq C\ell^{2\alpha-r}\|f\|_{\alpha}\|g\|_{\alpha}$ and 
\item[(iv)] $\|f-f*\varphi_{\ell}\|_0\leq C\ell\|f\|_1$.
\end{itemize}
The constants might depend on $r$ or $s$. Notice $\alpha=1$ the endpoint case is included.
\end{lemma}

\textsc{Remark:} if $f$ is not defined in all space then the domain where we can obtain the bound {\em shrinks an amount of $\ell$ units}. Notice also that for $s=0$ the constant in (i) can be taken to be exactly one.

We recall here a technical property that will be needed. 

\begin{lemma}\label{H4}
Given $G(u,v)$ a smooth function. Let $u,u',v,v'\in C^r$ with $\|u'\|_0+\|v'\|_0\leq C$ for a fixed constant. Then
\[\|G(u,v)-G(u',v')\|_r\leq C_r(\|u-u'\|_r+\|v-v'\|_r+(\|u'\|_r+\|v'\|_r)(\|u-u'\|_0+\|v-v'\|_0))\]
\end{lemma}
\textsc{Remark:} this appeared in K\"all\'en's work (cf. H4 condition in \cite{Ka}). We refer the reader there for a proof.

\section{Acknowledgments}

The author would like express his gratitude to C. De Lellis for drawing his attention to Kr\"oner's master thesis, his encouragement and enlightening discussions around the isometric embedding problem. 

This material is based upon work supported by the National Science Foundation under Grant No. DMS-1638352.

\end{document}